\documentclass[11pt]{article}
\usepackage[english]{babel}
\usepackage{amsmath}
\usepackage{amsthm}
\usepackage{amsfonts}
\usepackage{amssymb}
\usepackage{enumerate}
\usepackage{xspace}
\usepackage{euscript}
\usepackage{graphicx}
\usepackage{amscd}
\usepackage{epsfig}
\usepackage{tabularx}
\usepackage{url}
\usepackage{multimedia}

\begin{document}

\title{Recursive Integral Method with Cayley Transformation}

\author{R.~Huang, J.~Sun \\Department of Mathematical Sciences, Michigan Technological University \\ 
	C.~Yang \\ Computational Research Division, Lawrence Berkeley National Laboratory, Berkeley}


\maketitle

\begin{abstract}
Recently, a non-classical eigenvalue solver, called {\bf RIM}, was proposed to compute (all) eigenvalues in a region on the complex plane.
Without solving any eigenvalue problem, it tests if a region contains eigenvalues using an approximate spectral projection. 
Regions that contain eigenvalues are subdivided and tested recursively until eigenvalues are isolated with a specified precision. 
This makes {\bf RIM} an eigensolver distinct from all existing methods.
Furthermore, it requires no a priori spectral information.
In this paper, we propose an improved version of {\bf RIM} for non-Hermitian eigenvalue problems. 
Using Cayley transformation and Arnoldi's method, the computation cost is reduced significantly.
Effectiveness and efficiency of the new method are demonstrated by numerical examples and compared with 'eigs' in Matlab.
\end{abstract}
\section{Introduction}
We consider the non-Hermitian eigenvalue problem 
\begin{equation} \label{AxLambdaBx}
A x= \lambda B x,
\end{equation}
where $A$ and $B$ are $n\times n$ large sparse matrices. Here $B$ can be singular. Such eigenvalue problems arise in many scientific and engineering
applications \cite{GolubVorst2000CAM, Saad2011, SunZhou2016} as well as in emerging areas such as data analysis in social networks \cite{BoccalettiEtal2006PR}.

The problem of interest in this paper is to find (all) eigenvalues in a given region $S$ on the complex plane $\mathbb C$ without any spectral information,
i.e., the number and distribution of eigenvalues in $S$ are not known. 

In a recent work \cite{Huang2016JCP}, we developed an eigenvalue solver {\bf RIM} (recursive integral method).
\textbf{RIM}, which is essentially different from all the existing eigensolvers, is based on spectral projection and domain decomposition. Briefly speaking, 
given a region $S \subset \mathbb C$ whose boundary $\Gamma:=\partial S$ is a simple closed curve, 
{\bf RIM} computes an indicator $\delta_S$ using spectral projection $P$ defined by a Cauchy contour integral on $\Gamma$.
The indicator is used to decide if $S$ contains eigenvalue(s). When the answer is positive, $S$ is divided into sub-regions
and indicators for these sub-regions are computed. The procedure continues until the size of the region is smaller than 
the specified precision $\epsilon$ (e.g., $\epsilon = 10^{-6}$).
The centers of the regions are the approximations of eigenvalues.

To be specific, for $z \in \mathbb C$, the resolvent of the matrix pencil $(A,B)$ is defined as as
\begin{equation}
R_z(A,B):=(A-zB)^{-1}.
\end{equation}
Let $\Gamma$ be a simple closed curve lying in the resolvent set of $(A,B)$ on $\mathbb C$. Spectral projection for \eqref{AxLambdaBx} is given by
\begin{equation}
P(A,B)=\dfrac{1}{2\pi i}\int_{\Gamma}(A-zB)^{-1}dz.
\end{equation}
Given a random vector ${\boldsymbol f}$, it is well-known that $P$ projects ${\boldsymbol f}$
onto the generalized eigenspace associated with the eigenvalues enclosed by $\Gamma$. Clearly, $P{\boldsymbol f}$ is zero if there is no eigenvalue(s)
inside $S$, and nonzero otherwise.


{\bf RIM} differs from classical eigensolvers \cite{GolubVorst2000CAM} and recently developed integral based methods\cite{SakuraiSugiura2003CAM,Polizzi2009PRB}.
It simply computes an indictor of a region using the approximation to $P{\boldsymbol f}$,
\begin{equation}\label{XLXf}
P{\boldsymbol f} \approx  \dfrac{1}{2 \pi i} \sum_{j=1}^W \omega_j {\boldsymbol x}_j,
\end{equation}
where $\omega_j$'s are quadrature weights and ${\boldsymbol x}_j$'s are the solutions of the linear systems
\begin{equation}\label{linearsys}
(A- z_jB){\boldsymbol x}_j = {\boldsymbol f}, \quad j = 1, \ldots, W.
\end{equation}
Recall that if there is no eigenvalue inside $\Gamma$, then 
$P{\boldsymbol f} = {\bf 0}$ for all 
${\boldsymbol f} \in \mathbb C^n$.

In practice, one needs a threshold to distinguish 
between $|P{\boldsymbol f}| \ne 0$ and 
$|P{\boldsymbol f}|= 0$. The indicator needs to be robust enough to treat the following problems:
\begin{itemize}
\item[P1)] The randomly selected ${\boldsymbol f}$ may have only a 
small component in the range of $P$  (which we write as 
${\mathcal R}(P)$),
in which case $|P{\boldsymbol f}|$ may be small even 
when there are eigenvalues in $S$; 
\item[P2)] Since a quadrature rule is
used to approximate  $|P{\boldsymbol f}|$,
the indicator will not be zero (and may not 
even be very small) when there is no eigenvalue in $S$.
\end{itemize}
The strategy of {\bf RIM} in \cite{Huang2016JCP} selects a small threshold $\delta_0=0.1$ 
based on substantial experimentation, i.e., $S$ contains no eigenvalue if $\delta_S < \delta_0$. 
This choice of threshold 
for discarding a region systematically leans towards further 
investigation of regions that may potentially contain eigenvalues.
Of course, such a strategy leads to more computation cost.
 
To understand the first problem, consider an orthonormal basis 
$\{{\boldsymbol \phi}_j, j=1,\ldots, M\}$
for  ${\mathcal R}(P)$, which coincides with the
eigenspace associated with all the eigenvalues within $S$.
For a random ${\boldsymbol f}$,
\begin{equation}\label{fRP}
|P{\boldsymbol f}|=\left|{\boldsymbol f}|_{{\mathcal R}(P)}\right | =\left| \sum_{j=1}^M a_j {\boldsymbol \phi}_j \right |=\left(\sum_{j=1}^M a^2_i \right)^{1/2},
\end{equation}
where $a_j = ({\boldsymbol f}, {\boldsymbol \phi}_j)$. Since ${\boldsymbol f}$ is random, $|P{\boldsymbol f}|$ could be very small.
The solution in \cite{Huang2016JCP} is to normalize
$P{\boldsymbol f}$ and project again. 
The indicator $\delta_S$ is set to be
\begin{equation}
\label{sigmaP2P}
\delta_S := \left | P \left( \frac{P {\boldsymbol f}}{|P {\boldsymbol f}|}\right)\right|.
\end{equation}
Analytically, $P^2{\boldsymbol f} = P{\boldsymbol f}$. But numerical approximations to $P^2{\boldsymbol f}$  and $P{\boldsymbol f}$ may
differ significantly.  
\vskip 0.1cm

The following is the basic algorithm for {\bf RIM}:
\begin{itemize}
\item[] {\bf RIM}$(A, B, S, \epsilon, \delta_0, {\boldsymbol f})$
\item[]{\bf Input:}  matrices $A, B$, region $S$, precision $\epsilon$, threshold $\delta_0$, random vector ${\boldsymbol f}$.
\item[]{\bf Output:}  generalized eigenvalue(s) $\lambda$ inside $S$ 
\vskip 0.1cm
\item[1.] Compute ${\delta_S}$ using \eqref{sigmaP2P}.
\item[2.] Decide if $S$ contains eigenvalue(s).
	\begin{itemize}
		\item If $\delta_S < \delta_0$, then exit.
		\item Otherwise, compute the size $h(S)$ of $S$.
			\begin{itemize}
				\item[-] If $h(S)  > \epsilon $, 
						\begin{itemize}
						\item[] partition $S$ into subregions $S_j, j=1, \ldots N$.
						\item[] for $j=1: N$
						\item[] $\qquad${\bf RIM}$(A, B, S_j, \epsilon, \delta_0, {\boldsymbol f})$.
						\item[] end
						\end{itemize}
				\item[-] If $h(S) \le \epsilon$, 
						\begin{itemize}
							\item[] set $\lambda$ to be the center of $S$.
							\item[] output $\lambda$ and exit.
						\end{itemize}
			\end{itemize}
	\end{itemize}
\end{itemize}

To compute $\delta_S$, one needs to solve many linear systems 
\begin{equation}\label{fAzjByj}
  (A-z_j B) {\boldsymbol x}_j = {\boldsymbol f} 
\end{equation} 
parameterized by $z_j$.
In \cite{Huang2016JCP}, the Matlab linear solver `\textbackslash' is used to solve \eqref{fAzjByj}. This is certainly not efficient!

In this paper, we propose a new version of {\bf RIM}, called {\bf RIM-C}, to improve the efficiency.
The contributions include: 1) Cayley transformation and Arnoldi's method to speedup linear solves for the parameterized system
\eqref{fAzjByj}; and 2) a new indicator to improve the robustness and efficiency.
The rest of the paper is arranged as follows. In Section 2, we present how to incorporate Cayley transformation and the Arnoldi's method
into {\bf RIM}. In Section 3, we introduce a new indicator to decide if a region contains eigenvalues. Section 4 contains the new algorithm
and some implementation details. Numerical examples are presented in Section 5. We end up the paper with some conclusions and future
works in Section 6.

\section{Cayley Transformation and Arnoldi's Method}

\subsection{Cayley transformation}
The computation cost of {\bf RIM} mainly comes from solving the linear systems \eqref{fAzjByj} to compute the spectral projection $P{\boldsymbol f}$. 
In particular, when the method zooms in around an eigenvalue, it needs to solve linear systems for many close $z_j$'s. 
This is done one by one in the first version of {\bf RIM} \cite{Huang2016JCP}.
It is clear that the computation cost will be greatly reduced if one can take the advantage of the parametrized linear systems of same structure.

Without loss of generality, we consider a family of linear systems
\begin{equation} \label{eq:4}
(A-zB) {\boldsymbol x} ={\boldsymbol f},
\end{equation}
where $z$ is a complex number. 
When $B$ is nonsingular, multiplication of $B^{-1}$ on both sides of \eqref{eq:4} leads to
\begin{equation} \label{eq:Precond}
(B^{-1}A-z I) {\boldsymbol x}=B^{-1} {\boldsymbol f}.
\end{equation}
Given a matrix $M$, a vector ${\boldsymbol b}$, and a non-negative integer $m$, the Krylov subspace is defined as
\begin{align} 
K_{m}(M; {\boldsymbol b}):=\text{span} \{{\boldsymbol b}, M {\boldsymbol b}, \ldots, M^{m-1}{\boldsymbol b} \}.
\end{align}
The shift-invariant property of Krylov subspaces says that
\begin{align} \label{eq:5}
K_{m}(a M+ b I; {\boldsymbol b})=K_{m}(M; {\boldsymbol b}),
\end{align}
where $a$ and $b$ are two scalars.
Thus the Krylov subspace of $B^{-1}A-z I$ is the same as $B^{-1}A$, which is independent of $z$.

The above derivation fails when $B$ is singular.
Fortunately, this can be fixed by Cayley transformation \cite{Meerbergen2003SIAM}. 
Assume that $\sigma$ is not a generalized eigenvalue and $\sigma \neq z$. 
Multiplying  both sides of \eqref{eq:4} with
\begin{equation}\label{newprecond}
(A- \sigma B)^{-1},
\end{equation}
one obtains that
\begin{eqnarray*}
\label{ABsigmaz}(A- \sigma B)^{-1}{\boldsymbol f}&=&(A- \sigma B)^{-1}(A-z B){\boldsymbol x} \\
\nonumber &=&(A- \sigma B)^{-1}(A-\sigma B+(\sigma -z)B) {\boldsymbol x} \label{eq:Cayley} \\
\nonumber &=&(I+(\sigma - z)(A- \sigma B)^{-1}B) {\boldsymbol x}.\\
\end{eqnarray*}
Let $M=(A- \sigma B)^{-1}B$ and 
${\boldsymbol b}=(A- \sigma B)^{-1}{\boldsymbol f}$. Then
\eqref{eq:4} becomes
\begin{align} \label{eq:6}
(I+(\sigma -z)M) {\boldsymbol x} = {\boldsymbol b}.
\end{align}
From \eqref{eq:5}, the Krylov subspace $(I+(\sigma -z)M)$ is the same as $K_{m}(M; {\boldsymbol b})$. 

\subsection{Analysis of the pre-conditioners}
\begin{figure}[h!]
\begin{center}
\begin{tabular}{cc}
\resizebox{0.45\textwidth}{!}{\includegraphics{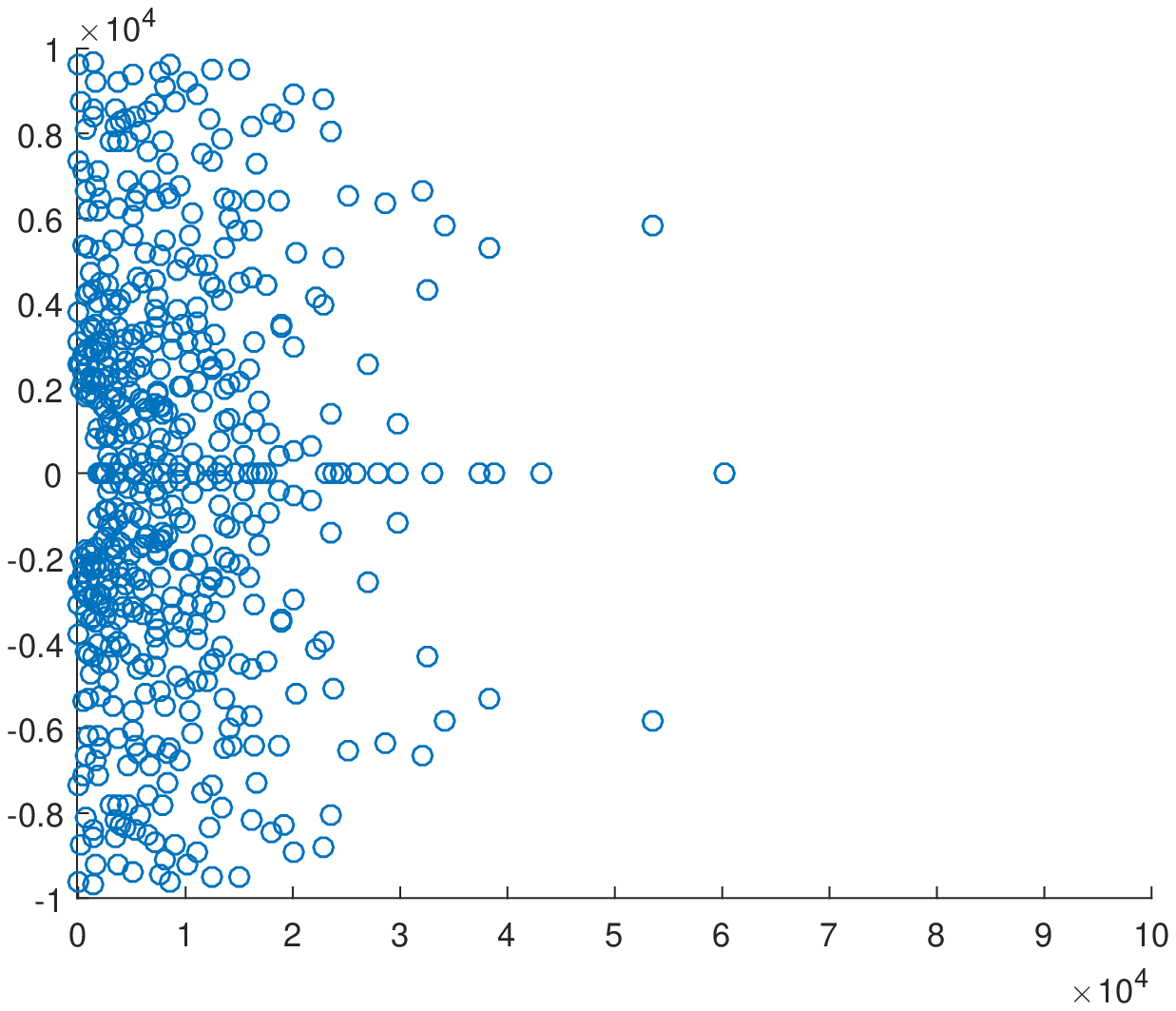}}&
\resizebox{0.45\textwidth}{!}{\includegraphics{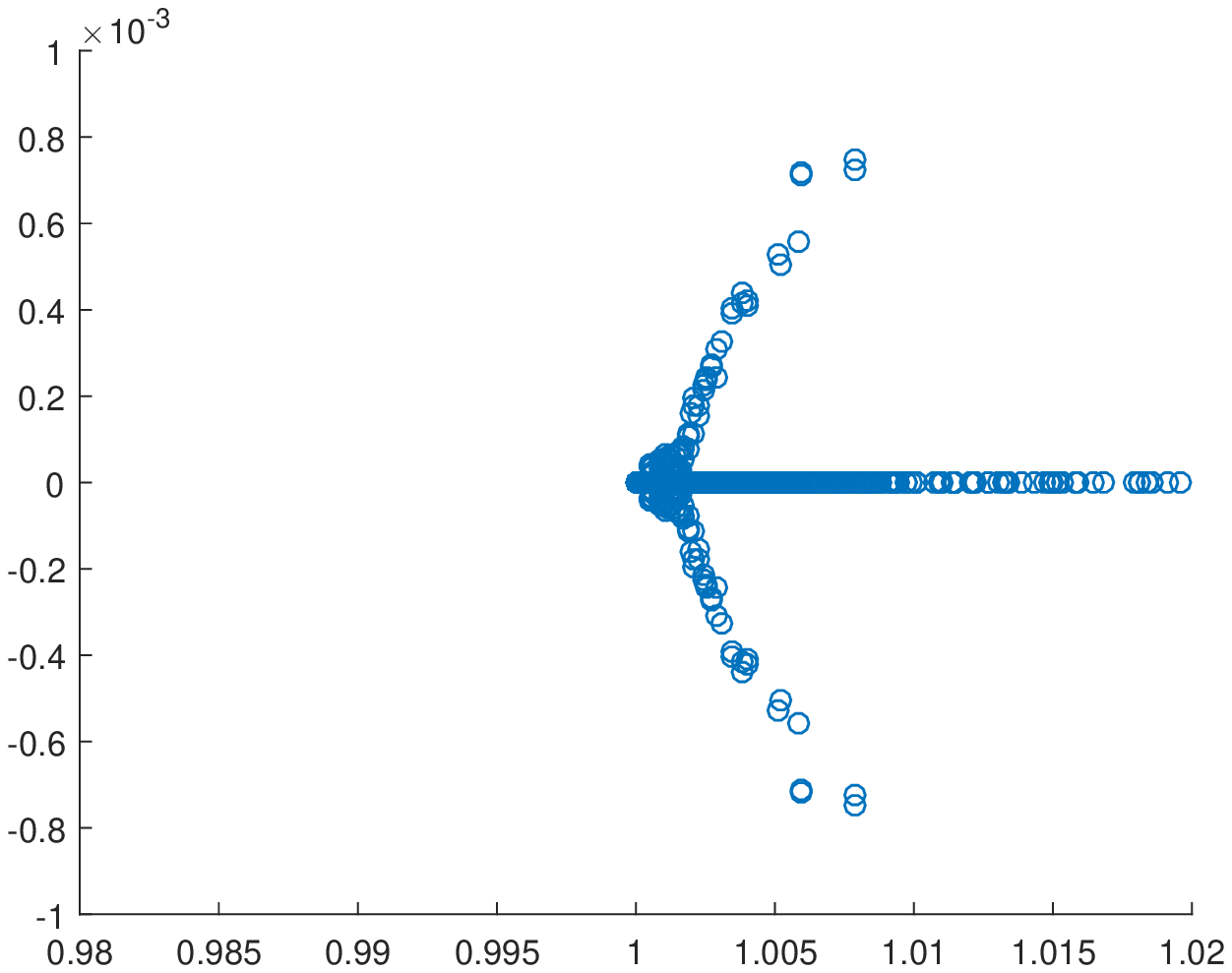}}
\end{tabular}
\end{center}
\caption{Matrices $A$ and $B$ are from Example 1 in Section 5. Left: Spectrum of original problem. Right: Spectrum after Cayley transformation.}
\label{SpectrumCT}
\end{figure}

Now we look at the connection between two pre-conditioners $B^{-1}$ and $(A- \sigma B)^{-1}$. 
Assume that $B$ is non-singular. Let $\lambda$ be an eigenvalue of $B^{-1}A$. Then $\theta=\dfrac{\lambda - z}{\lambda -\sigma}$ 
is an eigenvalue of 
\[
(A- \sigma B)^{-1}(A-z B).
\] 
The spectrum of $B^{-1}A$ might spread over the complex plane such that
Krylov subspace based iterative methods may not converge. However, after Cayley transformation, when $\lambda$ becomes large, 
$\theta$ will cluster around $1$ (see Fig.~\ref{SpectrumCT} for matrices $A$ and $B$ of {\bf Example 1} in Section 5).
Similar result holds when $B$ is singular. Note that when $\lambda$ approaches $\sigma$, $\theta$ will be very large in magnitude.
When $\lambda$ approaches $z$, $\theta$ goes to zero.
When $\lambda$ is away from $\sigma$ and $z$, $\theta$ is $O(1)$. The key here is that the spectrum of \eqref{eq:Cayley} 
has a cluster of eigenvalues around $1$ and only a few isolated eigenvalues, which favors fast convergence in Krylov subspace.

\subsection{Arnoldi Method for Linear Systems }
The computation cost can be significantly reduced by exploiting \eqref{eq:6}. Consider the orthogonal projection method for 
\[
M {\boldsymbol x}={\boldsymbol b}.
\]
Let the initial guess be ${\boldsymbol x}_0={\boldsymbol 0}$. One seeks an approximate solution ${\boldsymbol x}_m$ in 
$K_m(M; {\boldsymbol b})$ of dimension $m$ by imposing the Galerkin condition \cite{Saad2003}
\begin{equation} \label{eq:Krylov}
({\boldsymbol b}-M {\boldsymbol x}_m) \perp K_m(M; {\boldsymbol b}).
\end{equation}
The basic Arnoldi's process (Algorithm 6.1 of \cite{Saad2011}) is as follows.
\begin{itemize}
\item[1.] Choose a vector ${\boldsymbol v}_1$ of norm $1$
\item[2.] for $j=1,2, \ldots, m$
	\begin{itemize}
		\item $h_{ij} = (M{\boldsymbol v}_j, {\boldsymbol v}_i), \quad i=1,2, \ldots,j$,
		\item ${\boldsymbol w}_j = M {\boldsymbol v}_j - \sum_{i=1}^j h_{ij} {\boldsymbol v}_i$,
		\item $h_{j+1,j}=\|{\boldsymbol v}_j\|_2$, if $h_{j+1, j}=0$ stop
		\item ${\boldsymbol v}_{j+1} = {\boldsymbol w}_j/h_{j+1, j}$.
	\end{itemize}
\end{itemize}

Let $V_m$ be the $n \times m$ orthogonal matrix with column vectors ${\boldsymbol v}_1, \ldots, {\boldsymbol v}_m$ and $H_m$ be
the $m \times m$ Hessenberg matrix whose nonzero entries $h_{i,j}$ are defined as above.
From Proposition 6.6 of \cite{Saad2011}, one has that
\begin{equation} \label{eq:arnoldi}
M V_m=V_m H_m + {\boldsymbol v}_{m+1} {h}_{m+1,m}{\boldsymbol e}^{T}_m
\end{equation}
such that  
\[
\text{span}\{\mathrm{col}(V_m)\}=K_{m}(M; {\boldsymbol b}).
\] 
Let ${\boldsymbol x}_m=V_m {\boldsymbol y}$.
The Galerkin condition \eqref{eq:Krylov} becomes
\begin{equation}
V^{T}_m {\boldsymbol b}-V^{T}_m M V_m {\boldsymbol y} ={\boldsymbol 0}.
\end{equation}
Since $V_m^T MV_m = H_m$ (see Proposition 6.5 of \cite{Saad2003}), the following holds:
\[
H_m {\boldsymbol y}=V^{T}_m {\boldsymbol b}.
\]
From the construction of $V_m$, ${\boldsymbol v}_1=\dfrac{{\boldsymbol b}}{\|{\boldsymbol b}\|_2}$. Let $\beta=\|{\boldsymbol b}\|_2$. Then
\begin{equation}
      {\boldsymbol y}=\beta H_m^{-1}{\boldsymbol e}_1.
\end{equation}
Consequently, the residual of the approximated solution ${\boldsymbol x}_m$ can be written as
\begin{equation} \label{eq:err1}
\|{\boldsymbol b}-M {\boldsymbol x}_m\|_2 = {h}_{m+1,m}|{\boldsymbol e}_m^{T} {\boldsymbol y}|.
\end{equation}

Due to the shift invariant property, one has that
\begin{equation} \label{eq:arno}
\{I+(\sigma -z)M \} V_m =V_m( I+(\sigma -z)H_m)\\
+(\sigma - z){\boldsymbol v}_{m+1} {h}_{m+1,m}{\boldsymbol e}^{T}_m.
\end{equation}
By imposing a Galerkin condition similar to \eqref{eq:Krylov}, we have that
\begin{equation}
V^{T}_m {\boldsymbol b}-V^{T}_m \{I+(\sigma -z)M \} V_m {\boldsymbol y} =0,
\end{equation}
which implies
\begin{equation} \label{eq:many}
\{I+(\sigma -z)H_m \} {\boldsymbol y} =\beta {\boldsymbol e}_1.
\end{equation}
From \eqref{eq:err1}, one has that
\begin{equation}  \label{eq:error}
\|b-\{I+(\sigma - z)M \} {\boldsymbol x}_m\|_2 =(\sigma -z){h}_{m+1,m}|{\boldsymbol e}_m^{T} {\boldsymbol y}|.
\end{equation}

Matrix $M$ is an $n \times n$ matrix and $H_m $ is an $m \times m$ upper Hessenberg matrix such that $m \ll n$. 
Once $H_m$ and $V_m$ are constructed by Arnoldi's process, 
they can be used to solve \eqref{eq:many} for different $z$'s with residual given by \eqref{eq:error}. 
The residual can be monitored with a little extra cost. 

Next we explain how the Arnoldi's process is incorporated in {\bf RIM}.
To solve \eqref{fAzjByj} for quadrature points $z_j$'s, one chooses a proper shift $\sigma$. Following \eqref{eq:6}, one has that
\begin{align}
(I+(\sigma -z_j)M) {\boldsymbol x}_j={\boldsymbol b},
\end{align}
where $M=(A-\sigma B)^{-1} B$  and  ${\boldsymbol b}=(A- \sigma B )^{-1} {\boldsymbol f}$. 

From \eqref{eq:arno} and \eqref{eq:many}, 
\begin{eqnarray} \label{eq:y}
{\boldsymbol y}_j&=&\beta (I+ (\sigma-z_j) H_m)^{-1}{\boldsymbol e}_1,  \\
\nonumber {\boldsymbol x}_j &\approx& V_m {\boldsymbol y}_j, \\
 \label{eq:reduced}
   P{\boldsymbol f} &\approx&  \dfrac{1}{2 \pi i} \sum w_j V_m {\boldsymbol y}_j.
\end{eqnarray}
Hence the Krylov subspace for $M=(A-\sigma B)^{-1} B$ can be used to solve many linear systems associated with $z_j$'s close to $\sigma$.

\section{An Efficient Indicator}
Another critical problem of {\bf RIM} is to how to define the indicator $\delta_S$. As seen above, 
the indicator in \cite{Huang2016JCP} defined by \eqref{sigmaP2P} is to project a random vector twice.
One needs to solve linear systems with different right hand sides, i.e.,
${\boldsymbol f}$ and $P{\boldsymbol f}/|P{\boldsymbol f}|$.
Consequently, two Krylov subspaces, rather than one, are constructed for a single shift $\sigma$.

In this section, we propose a new indicator that avoids the construction of two Krylov subspaces.
The indicator stills needs to resolve the two problems (P1 and P2) in Section 1.   
The idea is to approximate $|P {\boldsymbol f}|$ with different sets of trapezoidal quadrature points
by taking the advantage of the Cayley transformation and Arnoldi's method discussed in the previous section.

Let $P{\boldsymbol  f}|_{n}$ be the approximation of $ P {\boldsymbol f} $ with $n$ quadrature points. 
It is well-known that trapezoidal quadratures 
of a periodic function converges exponentially \cite[Section 4.6.5]{davis1984methods}, i.e., 
\begin{align*}
 \left|P{\boldsymbol f}- P{\boldsymbol f}|_n\right| = O(e^{-C n}),
  \end{align*} 
where C is a constant depending on ${\boldsymbol f}$. The spectral projection satisfies
  \[ P {\boldsymbol f}|_{n}\begin{cases}
  \neq {\boldsymbol 0} & \text{if there are eigenvalues inside } S, \\
  \approx {\boldsymbol 0} & \text{no eigenvalue inside } S.
  \end{cases}  \]
For a large enough $n_0$, one has that
 \[ \dfrac{ \left | P {\boldsymbol f}|_{2n_0}\right|}{ \left | P {\boldsymbol f}|_{n_0} \right|}=\begin{cases}
  \dfrac{|P {\boldsymbol f}| + O(e^{-C 2n})}{|P {\boldsymbol f}| + O(e^{-C n})}  & \text{if there are eigenvalues inside } S, \\
 \dfrac{ O(e^{-C 2n})}{ O(e^{-C n})}=O(e^{-C n})  & \text{no eigenvalue inside } S.
  \end{cases}  \]
The new indicator is set to be
\begin{equation}\label{ISPf}
\delta_S = {|P {\boldsymbol f}_{2n_0}|}/{|P {\boldsymbol f}_{n_0}|}.
\end{equation}
A threshold value $\delta_0$ is also needed to decide if there exists eigenvalue in $S$ or not.
If $\delta_S > \delta_0:=0.2$, $S$ is said to be admissible, i.e., there exists eigenvalue(s) in $S$. The value $0.2$ is chosen based on numerical experimentation. 
Due to \eqref{eq:y} - \eqref{eq:reduced}, the computation cost to evaluate the new indicator is not expensive.

\section{The New Algorithm}
Now we are ready to give the algorithm in detail.
It starts with several shifts $\sigma$'s distributed in $S$ uniformly. The associated Krylov subspaces $K_m(M; {\boldsymbol b})$ are constructed and stored.
For a quadrature point $z$, the algorithm first attempts to solve the linear system \eqref{eq:4} using the Krylov subspace with shift $\sigma$ closest to $z$.
If the residual is larger than the given precision $\epsilon$, a Krylov subspace with a new shift $\sigma$ is constructed, stored and used to solve the
linear system.
Briefly speaking, the algorithm constructed some Krylov subspaces with different $\sigma$'s. 
These subspaces are then used to solve the linear system for all quadrature points $z_j$'s.  
From \eqref{eq:y} and \eqref{eq:reduced}, instead of solving a family of linear systems of size $n$, 
the algorithm solves linear systems of reduced size $m$ for most $z_j$'s.
This is the key idea to speed up {\bf RIM}.
We denote this improved version of {\bf RIM} by {\bf RIM-C} ({\bf RIM} with Cayley transformation).

Given a search region $S$ and a normalized random vector ${\boldsymbol f}$, we compute the indicator $\delta_S$ using \eqref{ISPf}.
Without loss of generality, $S$ is assumed to be a square. 
We set $n_0=4$ in \eqref{ISPf}.
If $\delta_S > 0.2$,  $S$ is divided uniformly into $4$ regions. The indicators 
of these regions are computed. This process continues until the size of the region is smaller than $d_0$.

\begin{enumerate}
\item[] {\bf Algorithm RIM-C:}
\item[] \textbf{RIM-C}$(A, B, S, {\boldsymbol f}, d_0, \epsilon, \delta_0, m, n_0)$
\item[] \textbf{Input:} 
	\begin{itemize}
		\item $A, B$: $n \times n$ matrices 
		\item $S$: search region in $\mathbb C$
		\item ${\boldsymbol f}$: a random vector
		\item $d_0$: precision 
		\item $\epsilon$: residual threshold 
		\item $\delta_0$: indicator threshold
		\item $m$: size of Krylov subspace
		\item $n_0$: number of quadrature points
	\end{itemize}
\item[] \textbf{Output:} 
	\begin{itemize}
	\item generalized eigenvalues inside $S$
	\end{itemize}
\end{enumerate}
\begin{enumerate}
\item Choose several $\sigma$'s uniformly in $S$ and construct Krylov subspaces 
\item Compute $\delta_S$ using \eqref{ISPf}. 
	\begin{itemize}
	\item[] Let $z$ be a quadrature point.
	\item Check if the linear system can be solved using the existing Krylov subspaces with residual less than $\epsilon$. 
	\item Otherwise, choose a new $\sigma$, construct a new Krylov subspace to solve the linear system.
	\end{itemize}
\item Decide if each $S$ contains eigenvalues(s).
\begin{itemize}
		\item If $\delta_S = \dfrac{|P {\boldsymbol f}|_{2n_0}|}{|P {\boldsymbol f}|_{n_0}|}<\delta_0$, exit.
		\item Compute the size of $S$, $h(S)$.			
				\item[-] If  $h(S)> \epsilon_0 $, uniformly partition $S_i$ into subregions $S_j, j=1, \ldots 4$
						\begin{itemize}
						\item[] for $j=1$ to $4$
						\item[] \qquad call \textbf{RIM-C}$(A, B, S_j, {\boldsymbol f}, d_0, \epsilon, \delta_0, m, n_0)$
						\item[] end
						\end{itemize}
				\item[-] Otherwise, output the eigenvalue $\lambda$ and exit.
			
	\end{itemize}
\end{enumerate}


\section{Numerical Examples}
In this section, {\bf RIM-C} (implemented in Matlab) is employed to compute all the eigenvalues in a given region.
To the authors' knowledge, there exists no eigensolver doing exactly the same thing. We compare {\bf RIM-C} with `eigs' in Matlab 
(IRAM: Implicitly Restarted Arnoldi Method \cite{arpack}). Although the comparison seems to be unfair to both methods, it gives some
idea about the performance of {\bf RIM-C}.

The matrices for {\bf Examples 1-5} come from a finite element discretization of the transmission eigenvalue problem \cite{JiSunTurner2012ACMTOM,Sun2011SIAMNA}
using different mesh size $h$. Therefore, the spectra of these problems are similar.
For Matlab function `eigs(A,B,K,SIGMA)', `K' and `SIGMA' denote the number of eigenvalues to compute and the {\it shift}, respectively.
For {\bf RIM-C}, the size of Krylov space is set to be $m=50$, $d_0 = 10^{-9}$, $\epsilon = 10^{-10}$, $\delta_0=0.2$, and $n_0=4$.
All the examples are computed on a Macbook pro with 16 Gb memory and 3 GHz Intel Core i7.
 
{\bf Example 1:}
The matrices $A$ and $B$ are $1018 \times 1018$ (mesh size $h\approx 0.1$). The search region $S=[1, 11] \times [-1, 1]$.
For `eigs', the `shift' is set to be $5.5$. For this problem, it is known that there exist $5$ eigenvalues in $S$. Therefore, `K' is set to be $5$.
Note that {\bf RIM-C} does not need this information.
The results are shown in Table~\ref{1018}.
Both {\bf RIM-C} and `eigs' compute $5$ eigenvalues and they are consistent. `eigs' uses less time than {\bf RIM-C}.

\begin{table}[h!]
\caption{Eigenvalues computed and CPU time by {\bf RIM-C} and `eigs' for {\bf Example 1}.}
\label{1018}
\centering
\begin{tabular}{l|r|r}
\hline
 & {\bf RIM-C} & `eigs' \\\hline
Eigenvalues & \textbf{3.9945390188}48445 & \textbf{3.9945390188}56096 \\
 & \textbf{6.93971914380}0903 & \textbf{6.93971914380}4773  \\
 & \textbf{6.9350539858}73570 & \textbf{6.9350539858}44678 \\ 
 & \textbf{10.6546658534}90588 & \textbf{10.6546658534}41946\\
 & \textbf{10.6587060246}50019 & \textbf{10.6587060246}09756\\
 \hline
CPU time& 0.284922s &0.247310s \\
\hline
\end{tabular}
\end{table}

{\bf Example 2:}
Matrices $A$ and $B$ are $4066 \times 4066$ (mesh size $h\approx 0.05$). Let $S=[20, 30] \times [-6, 6]$. For `eigs',  `{\it shift}' is set to be $25$.
Again, it is known in advance that there are $3$ eigenvalues in $S$. Hence `K' is set to be $3$.
The results are shown in Table~\ref{4066}. Both methods compute same eigenvalues and  `eigs' is faster.

\begin{table}[h!]
\caption{Eigenvalues computed and CPU time by {\bf RIM-C} and `eigs" for {\bf Example 2}.}
\label{4066}
\centering
\begin{tabular}{l|r|r}
\hline
 & {\bf RIM-C} & `eigs' \\\hline 
 Eigenvalues & \textbf{23.803023938}395199 $\qquad \qquad$ & \textbf{23.803023938}403236$\qquad \qquad$ \\
 & $\pm$  \textbf{5.6823043148}76092i & $\qquad \pm$  \textbf{5.6823043148}40053i\\ 
& \textbf{24.73702749700}6540 & \textbf{24.73702749700}3453\\
  & \textbf{24.7509596350}36583 & \textbf{24.7509596350}22376\\
 & \textbf{25.2781451874}65789 & \textbf{25.2781451874}57707\\
  & \textbf{25.2845015150}28143 & \textbf{25.2845015150}36474\\ 
 \hline
CPU time& 0.558687s &0.333513s \\
\hline
\end{tabular}
\end{table}

{\bf Example 3:}
Matrices $A$ and $B$ are $16258 \times 16258$ matrices (mesh size $h\approx 0.025$). 
Let $S=[0, 20] \times [-6, 6]$. There are $10$ eigenvalues in $S$. 
It is well-known that the performance of `eigs' is highly dependent on `shift'. 
In Table~\ref{16258},  we show the time used by {\bf RIM-C} and `eigs' with different shifts `{\it shift} =5, 10, 15'.
Notice that when the shift is not {\it good}, `eigs' uses much more time. In practice, {\it good} shifts are not known in advance.

\begin{table}[h!]
\caption{CPU time used by {\bf RIM} and `eigs' with different shifts for {\bf Example 3}.}
\label{16258}
\centering
\begin{tabular}{c|c|c|c|c}
\hline
 & {\bf RIM-C} & `eigs' shift=5& `eigs' shift=10 & `eigs' shift=15 \\ \hline
 CPU time& 2.571800s &0.590186 &\textbf{7.183679}s&0.392902s\\
 \hline
\end{tabular}
\end{table}


{\bf Example 4:} We consider a larger problem: $A$ and $B$ are  $260098 \times 260098$. 
Let $S=[0, 20] \times [-6, 6]$ (mesh size $h\approx 0.00625$). There are $17$ eigenvalues in $S$.
The results are in Table~\ref{260098}. This example, again, shows that for larger problems without any spectrum information,
the performance of {\bf RIM-C} is quite stable and consistent. However, the performance of `eigs' varies a lot with different `{\it shifts}'.
\begin{table}[h!]
\caption{CPU time used by {\bf RIM} and `eigs' with different {\it shifts} for {\bf Example 4}.}
\label{260098}
\centering
\begin{tabular}{c|c|c|c}
\hline
 & RIM-C & `eigs' shift=5 &`eigs' shift=10 \\ \hline 
CPU time& 104.228413s &\textbf{1696.703477}s&{\bf 272.506573}s\\
\hline
\end{tabular}
\end{table}

{\bf Example 5:} This example demonstrates the effectiveness and robustness of the new indicator.
The same matrices in {\bf Example 3} ($16258 \times 16258$) are used. Consider three regions $S_1, S_2$ and $S_3$.
$S_1=[18.4,18.8] \times [-0.2,0.2]$ has one eigenvalue inside. 
$S_2=[14.6,14.8] \times [-0.1,0.1]$ has two eigenvalues inside.
$S_3=[19.7,19.9] \times [-0.1,0.1]$ contains no eigenvalue.
Table~\ref{NewIndictorsV} shows the indicators of these three regions computed using \eqref{ISPf}.
It is seen that the indicator is different when there are eigenvalues inside the region and when there are no eigenvalues.
\begin{table}[h!]
\caption{Indictors: $S_1$ and $S_2$ contain at least one eigenvalue, $S_3$ contains no eigenvalue.}
\label{NewIndictorsV}
\centering
\begin{tabular}{c|r|r|r}
\hline
number of quadrature points & $P{\boldsymbol f}_{S_1}$ & $P{\boldsymbol f}_{S_2}$ &$P{\boldsymbol f}_{S_3}$\\
\hline
4 & 0.021036161440 & 0.000256531878& 0.001173702609\\
8 & 0.020981705584 & 0.000258504259 & 0.000044238403\\
\hline
$\delta_S$ & 0.997411 &0.992370 & 0.037691\\
\hline
\end{tabular}
\end{table}


Table~\ref{randv} shows the means, minima, maxima, and standard deviations of indicators of these three regions computed using $100$ random vectors.
The indicators are consistent for different random vectors.
\begin{table}[h!]
\caption{Means, minima, maxima, and standard deviations of indicators using $100$ random vectors.}
\label{randv}
\begin{center}
\begin{tabular}{lrrrr}
\hline
$S$ & mean  & min.&max.&std. dev.\\
\hline
 $S_1$&0.99848393687 &0.66250246918&1.43123449889&0.08740952445\\
$S_2$ &0.99926772105 &0.92600650392&1.14648387384&0.01788832121\\
$S_3$ &0.03763601782 &0.03734608324&0.03775912970&0.00010228556\\
\hline
\end{tabular}
\end{center}
\end{table}

\begin{figure}
\begin{center}
{ \scalebox{0.5} {\includegraphics{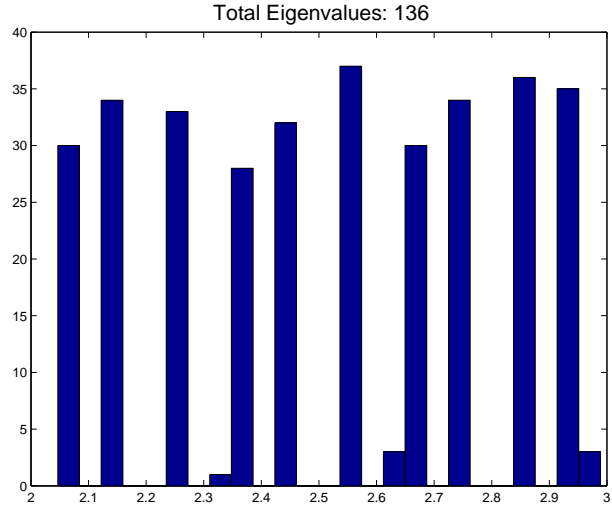}}}
\caption{Distribution of eigenvalues in $(2,3)$ for {\bf Example 6}.}
\label{Arnold}
\end{center}
\end{figure}

{\bf Example 6:} The last example shows the potential of {\bf RIM-C} to treat large matrices.
The sparse matrices are of ${\bf15,728,640 \times 15,728,640}$ arising from a finite element discretization of localized quantum states in 
random media \cite{Arnold2016}. {\bf RIM-C} computed $136$ real eigenvalues in $(2, 3)$, shown in the right picture
of Fig.~\ref{Arnold}. 

\section{Conclusions and Future Works}
This purposes of this paper is to compute (all) the eigenvalues of a large sparse non-Hermitian problem
in a given region. We propose a new eigensolver {\bf RIM-C}, which is an improved version of the recursive
integral method using spectrum projection. {\bf RIM-C} uses Cayley transformation and Arnoldi method
to reduce the computation cost.

To the authors' knowledge, {\bf RIM-C} is the only eigensolver for this particular purpose. As we mentioned,
the comparison of {\bf RIM-C} and `eigs' is unfair to both methods. However, the numerical results do show that {\bf RIM-C} 
is effective and has the potential to treat large scale problems.

Currently, the algorithm is implemented in Matlab. A parallel version using C++ is under development. 
For the time being, {\bf RIM-C} only computes eigenvalues, which is good enough for some applications. However, adding
a component to give the associated eigenvectors is necessary for other applications. It would also be useful to provide the
multiplicity of an eigenvalue. These are our future works to make {\bf RIM-C} a robust efficient eigensolver.

\end{document}